
%
\input amstex
\documentstyle{amsppt}

%
\define\[{\left[}
\define\]{\right]}
\define\({\left(}
\define\){\right)}
\define\<{\left<}
\define\>{\right>}
\define\QED{\ifmmode\qed\else{{\unskip\nobreak\hfil
\penalty50\null\nobreak\hfil$\qed$
\parfillskip=0pt\finalhyphendemerits=0\endgraf}}\fi}
\define\R{\widehat R}
\define\a{\alpha}
\redefine\b{\beta}
\define\g{\gamma}
\define\e{\varepsilon}
\define\s#1{\tilde#1}
\redefine\ss#1{\Tilde{\Tilde#1}}
\define\Diff{\operatorname{Diff}}

\define\coker{\operatorname{coker}}
\define\End{\operatorname{End}}
\define\Hom{\operatorname{Hom}}
\define\ch{\operatorname{char}}
%
%
\mag=\magstep 1
\pagewidth{15 true cm}
\pageheight{23 true cm}
\hcorrection{0.6 true cm}
\vcorrection{0.9 true cm}
\addto\tenpoint{\normalbaselineskip=1.5\normalbaselineskip
\normalbaselines}
\addto\eightpoint{\normalbaselineskip=1.5\normalbaselineskip
\normalbaselines}
%
%
\nopagenumbers
\pageno=1
\TagsOnRight
\chardef\tempcat=\the\catcode`\@\catcode`\@=11
  \def\logo@{}
\catcode`\@=\tempcat
\refstyle{C}
\widestnumber\no{99}
%
%

%
%
\topmatter
\title On quantized algebra of Wess-Zumino differential
operators at roots of unity
\endtitle
\author  Alexander Verbovetsky \endauthor
\rightheadtext{Quantized differential operators}
\affil
Scuola Internazionale Superiore di Studi Avanzati\\
Trieste, Italy
\endaffil
\date {3 May 1995} \enddate
\address
S.I.S.S.A., Via Beirut 2-4, 34013 Trieste, Italy
\endaddress
\email verbovet\@sissa.it \endemail
\endtopmatter
\newpage
\document
\subhead 1. Motivation \endsubhead
In \cite{13} J.~Wess and B.~Zumino constructed a broad family of
quantum deformations of the algebra of scalar differential operators on
the affine space.  If the deformation parameter $q$ is generic, this
algebra has the same dimension as in non-deformed case; if $q$ is a
root of unity, then the dimension of deformed algebra is less than that
of non-deformed algebra.  However, this does not quite agree with the
intuition gained in the classical theory of differential operators.  In
fact, due to the analogy between quantum case at roots of unity and the
classical case in positive characteristic, one should expect the
appearance of differential operators that are not compositions of
operators of order $\leqslant 1$, whereas the dimension of the whole
algebra must be the same as in non-deformed case.  In this short paper
it is shown that this intuition is a true one:  we quantize the
standard algebraic definition of differential operator
$\[\ldots\[\[D,f_0\],f_1\]\ldots,f_k\]=0$ by replacing the usual
commutators by twisted ones and obtain the algebra of differential
operators that has the classical dimension for any $q$ and coincides
with the Wess-Zumino algebra at generic parameter values.  A detailed
description of this algebra is presented.  Note also that having such
quantum differential operators one can construct quantum jets, de~Rham
and Spencer complexes, integral forms, Euler operator, and so on,
following the same logic as in non-deformed case (see \cite{9, 10, 12,
11, 7}).  A realization of this program has been started in \cite{8}.

\subhead 2 \endsubhead
Here we recall some standard facts about quantum affine spaces. We
refer to \cite{6, 2, 1, 4} for further details.

Let $\Bbbk$ be a commutative ground ring with unity.  Consider a free
$\Bbbk$-module $V$ and a non-degenerate linear operator $\R\:V\otimes
V\to V\otimes V$.  Throughout this paper we assume that the operator
$\R$ satisfies the Yang-Baxter equation
$$\R_{12}\R_{23}\R_{12}=\R_{23}\R_{12}\R_{23}$$
(here $\R_{12}=\R\otimes 1$, $\R_{23}=1\otimes\R)$ and the Hecke
condition
$$\R^2=\(q-q^{-1}\)\R+1$$
for some invertible element $q$ of $\Bbbk$.

\example{Basic example (\cite{3})}
Let $\{\xi_1,\ldots,\xi_n\}$ be a basis of $V$ and $\{x^1,\ldots,x^n\}$
be the dual basis of $V^*$. Let $\|\R^{ij}_{kl}\|$ be the matrix of
$\R$ in this basis. Choose an invertible $q\in\Bbbk$ and put
$$\R^{ij}_{kl}=\delta^i_l\delta^j_k\(1+\(q-1\)\delta^{ij}\)+
\(q-q^{-1}\)\delta^i_k\delta^j_l\theta\(j-i\), \tag{1}$$
where $\theta\(i\)=\cases1,&\text{for $i>0$}\\0,&\text{for
$i\leqslant0$}. \endcases$
\endexample

The algebra $A$ of a {\it quantum affine space\/} is defined as the
quotient algebra of the tensor algebra $T\(V^*\)$ by the ideal
generated by the image of $q-\R^*\:V^*\otimes V^*\to V^*\otimes V^*$.
In coordinates, $A$ is generated by $x^1,\ldots,x^n$ subject to the
relations
$$\R^{ij}_{\a\b}x^{\a}x^{\b}=qx^ix^j$$
(we assume summation over repeated indices which occur in both upper
and lower positions).

\example{Example}
Assume that $\R$ is given by \thetag{1}. Then $A$ is generated by
$x^1,\ldots,x^n$ modulo the following relations:
$$x^ix^j=qx^jx^i,\quad i<j.$$
\endexample

Consider the matrix algebra $M$ defined as the quotient of the tensor
algebra $T\(\End\(V\otimes V\)\)$ by the ideal generated by the
elements $\R F-F\R$, where $F\in\End\(V\otimes V\)$. The algebra $M$ is
generated by $t^i_j$, $1\leqslant i,j\leqslant n$, obeying the
relations
$$\R^{ij}_{\a\b}t^{\a}_kt^{\b}_l=t^i_{\a}t^j_{\b}\R^{\a\b}_{kl}.$$
There is a natural bialgebra structure on $M$, with the
comultiplication given by $\Delta t^i_j=t^i_{\a}\otimes t^{\a}_j$ and
the counit $\e t^i_j=\delta^i_j$. One obtains a matrix {\it quantum
group} $H$ by taking the Hopf envelop of $M$, that is an initial
object in a category of bialgebra morphisms $M\to C$, with $C$ being
a Hopf algebra (see \cite{5}).

The algebra $A$ is a (right) comodule-algebra over $H$, with the
coaction given by
$$\Delta\(x^i\)=x^{\a}\otimes\s{t^i_{\a}},$$
where the tilde stands for the antipode.

\remark{Remark}The antipode was inserted here to make this
transformation {\it right\/} coaction.
\endremark

We shall assume that the Hopf envelop $H$ of the matrix bialgebra is
endowed with a cobraided structure, i\.e., that there
exists a non-degenerate bilinear form $\<\cdot\,,\cdot\>$
on $H$ satisfying the conditions:
$$\gather\aligned \<h,h_1h_2\>&=\<h_i,h_2\>\<h^i,h_1\> \\
                  \<h_1h_2,h\>&=\<h_1,h_i\>\<h_2,h^i\> \endaligned \\
              g_jh_i\<h^i,g^j\>=\<h_i,g_j\>h^ig^j \endgather$$
for all $h$, $h_1,h_2,g\in H$ and $\Delta h=h_i\otimes h^i$, $\Delta
g=g_j\otimes g^j$. On the generators this form will be given by
$$\<t^i_j,t^k_l\>=q(\R^{-1})^{ki}_{jl}.$$
It is well-known that such a structure does exist in all examples
one is likely to encounter, provided $H$ was constructed in an
appropriate category.

Having this form we can introduce a (left) $H$-module structure on the
algebra~$A$:
$$h\cdot f=f_i\<h,h^i\>,$$
where $h\in H$, $f\in A$, $\Delta f=f_i\otimes h^i$. In coordinates we
have
$$\align
t^i_jx^k&=C^{ik}_{\a j}x^{\a},\quad t^i_j\cdot 1=\delta^i_j, \\
\s{t^i_j}x^k&=q(\R^{-1})^{ki}_{j\a}x^{\a}, \\
\ss{t}^i_jx^k&=q^{-1}\R^{ik}_{\a j}x^{\a}, \endalign$$
where the matrix $C$ is defined by
$C^{i\a}_{j\b}(\R^{-1})^{k\b}_{l\a}=q^{-1}\delta^i_l\delta^k_j$. Note
that this makes $A$ a crossed $H$-bimodule.
\proclaim{Lemma}{\it
For any $f\in A$ we have
$$\align x^if&=\(\ss{t}^i_{\a}f\)x^{\a} \\
fx^i&=x^{\a}\(\s{t}^i_{\a}f\).\endalign$$}
\endproclaim
\demo{Proof}
We prove the first relation, the second can be proved in the same way.
For $f=1$ the relation is obvious. Suppose inductively that it is true
for all $f$ such that $\deg f<m$. Take $f=x^jg$ with $\deg g<m$. Then
we get
$$\(\ss{t}^i_{\a}f\)x^{\a}=\(\ss{t}^i_{\b}x^j\)\(\ss{t}^{\b}_{\a}g\)
x^{\a}=\(\ss{t}^i_{\b}x^j\)x^{\b}g=q^{-1}\R^{ij}_{\a\b}x^{\a}x^{\b}g
=x^ix^jg=x^if.$$
Therefore the relation holds true for any$f$ of degree $m$. \QED
\enddemo

Having an $H$-module structure on $A$ we can define an $H$-module
structure on $\Hom_{\Bbbk}\(A,A\)$ by the formula:
$$\(t^i_jD\)\(f\)=t^i_{\a}\(D\(\s{t}^{\a}_jf\)\),\quad D\in
\Hom_{\Bbbk}\(A,A\).$$
Notice that $\Hom_{\Bbbk}\(A,A\)$ is a module-algebra over $H$:
$$\(t^i_j\(D_1\circ D_2\)\)=t^i_{\a}D_1D_2\s{t}^{\a}_j=
t^i_{\b}D_1\s{t}^{\b}_{\e}t^{\e}_{\a}D_2t^{\a}_j=\(t^i_{\e}D_1\)\circ
\(t^{\e}_jD_2\).$$

\subhead 3 \endsubhead
Now we define quantum algebra of (scalar) differential operators. For
any $\Bbbk$-linear mapping $D\:A\to A$ and for all $k\geqslant0$,
$1\leqslant i\leqslant n$, define the mappings $\[D,x^i\]_k\:A\to A$ by
the formula:
$$\[D,x^i\]_k\(f\)=D\(x^if\)-q^{2k}x^{\a}\(\s{t}^i_{\a}D\)\(f\).$$
The set of differential operators $\Diff_k\(A\)$ of
order $\leqslant k$ can be defined by the following inductive
procedure:
$\Diff_0\(A\)=\{D\in\Hom_{\Bbbk}\(A,A\)|\[D,x^i\]_0\allowmathbreak=0
\medspace\forall i\}$. In general, $\Diff_{k+1}\(A\)=\{D\in
\Hom_{\Bbbk}\(A,A\)|\[D,x^i\]_{k+1}\in\Diff_k\(A\)\medspace\forall
i\}$.
\proclaim{Proposition}{\it \roster
\runinitem If $D\in\Diff_k\(A\)$ then $h\cdot
D\in\Diff_k\(A\)\quad\forall h\in H$.
\item If $D_1\in\Diff_k\(A\)$ and $D_2\in\Diff_l\(A\)$, then
$D_1\circ D_2\in\Diff_{k+l}\(A\)$.
\endroster}\endproclaim
\demo{Proof}\therosteritem1 It is sufficient to prove that
$$t^j_l\[D,x^i\]_k=C^{\a i}_{\b l}\[t^j_{\a}D,x^{\b}\]_k.$$
We have
$$\split
t^j_l\[D,x^i\]_k&=\(t^j_{\a}D\)\(t^{\a}_lx^i\)-q^{2k}\(t^j_{\b}x^{\a}
\)\(t^{\b}_l\s{t}^i_{\a}D\)\\
&=C^{\a i}_{\b l}\(t^j_{\a}D\)x^{\b}-q^{2k}C^{j\a}_{\e\b}
x^{\e}\(t^{\b}_l\s{t}^i_{\a}D\)\\
&=C^{\a i}_{\b l}\(t^j_{\a}D\)x^{\b}-q^{2k}C^{\a i}_{\b
l}x^{\e}\(\s{t}^{\b}_{\e}t^j_{\a}D\)\\
&=C^{\a i}_{\b l}\[t^j_{\a}D,x^{\b}\]_k.\endsplit$$
$$\multline\text{\therosteritem2 }\[D_1\circ D_2,x^i\]_{k+l}\\
=D_1D_2x^i-q^{2l}D_1x^{\a}\(\s{t}^i_{\a}D_2\)+
q^{2l}D_1x^{\a}\(\s{t}^i_{\a}D_2\)-q^{2k+2l}x^{\b}\(\s{t}^{\a}_{\b}D_1
\)\(\s{t}^i_{\a}D_2\)\\
=D_1\circ\[D_2,x^i\]_l+q^{2l}\[D_1,x^{\a}\]_k\circ\s{t}^i_{\a}D_2
\endmultline$$
This implies the statement in a standard way. \QED
\enddemo
This Proposition has the following:
\proclaim{Corollary}{\it \roster
\runinitem There exists an $A$-bimodule structure on $\Diff_k\(A\)$
\rom{:}
$$a\cdot D=a\circ D,\quad D\cdot a=D\circ a.$$
\item For any $k<l$ we have $\Diff_k\(A\)\subset\Diff_l\(A\)$.
\endroster}\endproclaim

We denote $\Diff\(A\)=\bigcup_{k\geqslant0}\Diff_k\(A\)$.

\subhead 4 \endsubhead
Now we give a detailed description of the modules $\Diff_k\(A\)$. We
shall assume that the space $V^*$ has an ordered basis $x^1,\ldots,x^n$
such that the monomials $\(x^1\)^{i_1}\ldots\(x^n\)^{i_n}$ form a basis
of $A$ and $\R^{ij}_{kl}=0$ if $i>j$ and either $k>l$ or $k\geqslant
i$.

Clearly, any differential operator $D\in\Diff_k\(A\)$ is completely
determined by its values on the basic monomials of degree $\leqslant
k$. The following theorem states that these values can be taken in an
arbitrary way.

\proclaim{Theorem}{\it The $A$-module $\Diff_k\(A\)$ has a basis $1$,
$\partial_{\sigma}$, $|\sigma|\leqslant k$, where $\sigma=j_1\ldots
j_r$, $1\leqslant j_l\leqslant j_{l+1}\leqslant n$, is a multi-index,
$|\sigma|=r$, and the operator $\partial_{\sigma}$ is equal to $1$ at
$x^{j_1}\ldots x^{j_r}$ and vanishes on the other monomials of degree
$\leqslant k$.}
\endproclaim

\demo{Proof} We have to show that the defining relations $\[\cdots\[D,
x^{i_k}\]_k\cdots x^{i_0}\]_0=0$ with $i_{l+1}\leqslant i_l$ will be
sufficient for $D\in\Diff_k\(A\)$. For this we verify that
$$\[\[D,x^i\]_{l+1},x^j\]_l=q^{-1}\R^{ij}_{\a\b}\[\[D,x^{\a}\]_{l+1},
x^{\b}\]_l.$$
We have
$$\multline\[\[D,x^i\]_{l+1},x^j\]_l\\
\aligned &=Dx^ix^j-q^{2l+2}x^{\a}\(
\s{t}^i_{\a}D\)x^j-q^{2l}x^{\a}\(\s{t}^{\g}_{\a}D\)\(\s{t}^j_{\g}x^i\)
+q^{4l+2}x^{\a}\(\s{t}^{\g}_{\a}x^{\e}\)\(\s{t}^j_{\g}\s{t}^i_{\e}D\)\\
&=Dx^ix^j-q^{2l+1}\(q\delta^i_{\g}\delta^j_{\e}+(\R^{-1})^{ij}_{\g\e}\)
\(x^{\a}\(\s{t}^{\g}_{\a}D\)x^{\e}\)+q^{4l+2}x^{\e}x^{\g}\(\s{t}^j_{\g}
\s{t}^i_{\e}D\)\\&=q^{-1}\R^{ij}_{\a\b}Dx^{\a}x^{\b}-q^{2l}
\R^{ij}_{\a\b}\(q\delta^{\a}_{\g}\delta^{\b}_{\e}+(\R^{-1})^{\a\b}_
{\g\e}\)\(x^{\varkappa}\(\s{t}^{\g}_{\varkappa}D\)x^{\e}\)\endaligned\\
+q^{4l+1}\R^{ij}_{\a\b}x^{\e}x^{\g}\(\s{t}^{\b}_{\g}\s{t}^{\a}_{\e}D\)=
q^{-1}\R^{ij}_{\a\b}\[\[D,x^{\a}\]_{l+1},x^{\b}\]_l.\QED
\endmultline$$
\enddemo

\subhead 5 \endsubhead
We can define a (right) $H$-comodule structure on $\Diff_k\(A\)$ by the
following property:
$$h\cdot D=D_i\<h,h^i\>,$$
where $h\in H$, $D\in\Diff_k\(A\)$, $\Delta D=D_i\otimes h^i$. In
particular, for the operators $\partial_i$ this yields:
$$\Delta\partial_i=\partial_{\a}\otimes t^{\a}_i.$$
This structure makes $\Diff\(A\)$ an $H$-comodule-algebra:
$$\multline h\cdot\(D_1\circ
D_2\)=\(h_iD_1\)\(h^iD_2\)=\(D_1\)_k\(D_2\)_l \<h_i,g^k_1\>
\<h^i,g^l_2\>\\=\(D_1\)_k\(D_2\)_l\<h,g^l_2g^k_1\>,\endmultline$$
so that
$$\Delta\(D_1\circ D_2\)=\(D_1\)_k\circ\(D_2\)_l\otimes g^l_2g^k_1,$$
where $\Delta h=h_i\otimes h^i$, $\Delta D_r=\(D_r\)_k\otimes g^k_r$,
$r=1,2$.

\proclaim{Proposition (Covariance of differential operators)}{\it
For any $D\in\Diff_k\(A\)$ and $f\in A$
$$\Delta\(D\(f\)\)=\Delta\(D\)\(\Delta\(f\)\),$$
where the action in right-hand side is given by $\(\nabla\otimes h_1\)
\(a\otimes h_2\)=\nabla\(a\)\otimes h_2h_1$, $\nabla\in\Diff_k\(A\)$,
$a\in A$, $h_1,h_2\in H$.}
\endproclaim

\demo{Proof} For $f$ of degree zero the statement is trivial. Suppose
inductively that it is true for $\deg f\leqslant m$ and take an element
$f$ of degree $m$. We get
$$\multline\Delta\(D\(x^if\)\)=\Delta\(\(Dx^i\)\(f\)\)=\Delta\(Dx^i\)\(
\Delta\(f\)\)=\Delta\(D\)\Delta\(x^i\)\Delta\(f\)\\=\Delta\(D\)\(\Delta
\(x^if\)\).\QED\endmultline$$
\enddemo

\subhead 6 \endsubhead
Consider the first order operators $\partial_i$. It follows easily from
the very definition that they satisfy the Wess-Zumino Leibnitz rule:
$$\partial_ix^j=\delta^j_i+q\R^{j\b}_{i\a}x^{\a}\partial_{\b}.$$
Further, one has the following commutation relations between
$\partial_i$:
$$ \R^{\b\a}_{ji}\partial_{\a}\partial_{\b}=q\partial_i\partial_j.$$
Indeed, it is obvious that the second order operator
$\R^{\b\a}_{ji}\partial_{\a}\partial_{\b}-q\partial_i\partial_j$
vanishes on the monomials of degree $0$ and $1$. For the degree $2$ we
have:
$$\multline
\(\R^{\b\a}_{ji}\partial_{\a}\partial_{\b}-q\partial_i\partial_j\)\(x^k
x^l\)=\R^{kl}_{ji}+q\R^{\b\a}_{ji}\R^{kl}_{\b\a}-q\delta^l_i\delta^k_j-
q^2\R^{kl}_{ji}\\=q(\R^2-\(q-q^{-1}\)\R-1)^{kl}_{ji}=0.
\endmultline$$

\definition{Definition (\cite{2})}
The {\it quantum Weyl algebra\/} $\Cal D$ is the algebra with $2n$
generators $x^1,\ldots,x^n,\partial_1,\ldots,\partial_n$ satisfying the
following commutation relations:
$$\gather
\R^{ij}_{\a\b}x^{\a}x^{\b}=qx^ix^j, \\
\R^{\b\a}_{ji}\partial_{\a}\partial_{\b}=q\partial_i\partial_j, \\
\partial_ix^j=\delta^j_i+q\R^{j\b}_{i\a}x^{\a}\partial_{\b}.
\endgather$$
\enddefinition

{}From the previous discussion it follows that there exists an algebra
morphism $w\:\Cal D\to\Diff\(A\)$ such that $w\(x^i\)=x^i$,
$w\(\partial_j\)=\partial_j$. The image of $w$ consists of differential
operators that are compositions of operators of order $\leqslant1$. One
can also show that $\coker w$ is exactly the zero Spencer cohomology of
the algebra $A$ (see \cite{12}). For  the matrix $\R$ of the form
\thetag1 from the above-stated Theorem and a result of E.~Demidov
\cite{2, Corollary 12.6} it follows the following:

\proclaim{Proposition}{\it
If $q^2$ is not a root of unity, or $q^2=1$ and $\ch\Bbbk=0$, then
$w\:\Cal D\to\Diff\(A\)$ is an isomorphism, and if $q^2$ is a primitive
root of unity of degree $\ell$, then the kernel of $w$ is generated by
$\partial^{\ell}_1,\ldots,\partial^{\ell}_n$.}
\endproclaim

Thus, in the root of unity case there exist operators that can not be
presented as a composition of operators of order $\leqslant1$.

\example{Example}
Consider the operators $D_j=\partial_{\undersetbrace
\text{$\ell$ times} \to{j\ldots j}}$. It is easy to see that
$\partial_i\(\(x^j\)^\ell\)\allowmathbreak=0\quad\forall i$. Hence
$D_j$ is not a composition of operators of order $\leqslant1$. In fact,
any operator is a composition of $D_j$ and operators of order
$\leqslant1$.
\endexample

\subsubhead Acknowledgments \endsubsubhead
The author is grateful to the Italian Ministero degli Afari Esteri for
the financial support and to the SISSA for the excellent working
conditions.

\enddocument